\documentclass[11pt,reqno]{amsart}
\usepackage[a4paper,margin=27mm,headheight=14pt]{geometry}
\usepackage{amsmath,amssymb,amsfonts,amsthm,mathtools,booktabs,mathrsfs}
\usepackage{microtype}
\usepackage{xcolor}
\usepackage[colorlinks=true,linkcolor=blue!45!black,citecolor=blue!45!black,urlcolor=blue!45!black]{hyperref}
\newif\ifshowchanges
\showchangesfalse
\newenvironment{revision}{%
  \ifshowchanges\color{red}%
    \hypersetup{linkcolor=red,citecolor=red,urlcolor=red}%
  \fi}{\par}
\newcommand{\rev}[1]{\ifshowchanges{\color{red}#1}\else#1\fi}
\hypersetup{pdftitle={Optimal Sample Complexity for Genus-Two Zeta Functions from Noisy Random Euler Products},pdfauthor={Taro Hayashi and Ryoichi Suzuki}}
\newtheorem{theorem}{Theorem}[section]
\newtheorem{proposition}[theorem]{Proposition}
\newtheorem{lemma}[theorem]{Lemma}
\newtheorem{corollary}[theorem]{Corollary}
\theoremstyle{definition}
\theoremstyle{remark}
\newcommand{\F}{\mathbb F}\newcommand{\R}{\mathbb R}
\newcommand{\E}{\mathbb E}\newcommand{\Var}{\operatorname{Var}}
\newcommand{\Prob}{\mathbb P}\newcommand{\Law}{\operatorname{Law}}
\newcommand{\TV}{\operatorname{TV}}\newcommand{\He}{\mathsf H}
\newcommand{\at}{\operatorname{arctanh}}\newcommand{\dd}{\,\mathrm d}
\newcommand{\eps}{\varepsilon}

\DeclareMathOperator{\Arg}{Arg}
\allowdisplaybreaks[2]
\numberwithin{equation}{section}

\title[Finite-sample zeta rigidity]{Optimal Sample Complexity for Genus-Two Zeta Functions from Noisy Random Euler Products}
\author{Taro Hayashi}
\address{Department of Mathematical Sciences, College of Science and Engineering, Ritsumeikan University, Kusatsu, Shiga 525-8577, Japan}
\email{haya4@fc.ritsumei.ac.jp}
\author{Ryoichi Suzuki}
\address{Department of Business Economics, School of Management, Tokyo University of Science, 1-11-2 Fujimi, Chiyoda, Tokyo 102-0071, Japan and Research Center for Finance,
BKC Research Organization of Social Sciences,
Ritsumeikan University,
Kusatsu, Shiga 525-8577, Japan.}
\email{rsuzukimath@gmail.com}
\thanks{ORCID IDs:
Taro Hayashi,
\href{https://orcid.org/0009-0005-0145-8042}
{0009-0005-0145-8042};
Ryoichi Suzuki,
\href{https://orcid.org/0000-0001-9979-1882}
{0000-0001-9979-1882}.}
\date{}
\begin{document}
\begin{abstract}
We study exact recovery of the Hasse--Weil zeta numerator from repeated noisy observations of a randomized Euler product. For a genus-two curve over $\F_q$, we attach independent symmetric signs to its closed points and observe the logarithm of the resulting product at $t_q=c/\sqrt q$, where $0<c<1/\sqrt{12}$ is fixed. For every fixed, known additive noise density $f$ satisfying $\sqrt f\in H^1(\R)$, we prove that the minimax number of independent observations required for exact recovery with error probability at most $\delta$ has order $q^4\log(1/\delta)$. No moment condition is imposed on the noise, so the result includes Cauchy readout. The upper bound estimates the phase of an empirical characteristic function and uses the separation of the integral Weil coefficients. The matching lower bound compares the full observation laws for two actual curves with ordinary, absolutely simple Jacobians. It combines a binomial comparison with a uniform bound on the closed-point tail. Thus access to the full observations does not improve the sample order attained by the phase estimator. The result concerns the information in a specified observation channel. It does not assume an efficient implementation of that channel.
\end{abstract}
\maketitle

\section{Introduction}
Let $C/\F_q$ be a smooth, projective, geometrically connected curve of genus two. Its Hasse--Weil zeta function is
\[
 \zeta_C(T)=\exp\left(\sum_{n\ge1}\frac{\#C(\F_{q^n})}{n}T^n\right)
 =\frac{P_C(T)}{(1-T)(1-qT)},
\]
where, with the sign convention used throughout,
\begin{equation}\label{eq:poly}
 P_C(T)=1+a_CT+b_CT^2+qa_CT^3+q^2T^4\in\mathbb Z[T].
\end{equation}
The classical Weil bounds give $|a_C|\le4\sqrt q$ and $|b_C|\le6q$. Determining $(a_C,b_C)$ is equivalent to determining the zeta function, but it does not determine the curve up to isomorphism.

A random Euler product compresses all closed-point counts into a random analytic quantity. Observing it at one real argument leaves a scalar inference problem: can the integral Frobenius data still be recovered, and does using the full observation improve on recovering a single location parameter? We answer these questions by matching a constructive recovery bound with a lower bound for every measurable estimator. The lower bound is essential here: a small difference of means alone does not imply that two observation laws are hard to distinguish.

Random Euler products also occur in models of function-field $L$-functions. Sawin~\cite{Sawin} combines a random Euler product with a random-matrix model to approximate statistics of a family of $L$-functions. Here the curve is the unknown parameter, and the random Euler product is an observation channel. Our sample count measures statistical information in that channel.

Write $|C|$ for the closed points of $C$, and let $\eps_x\in\{-1,1\}$ be mutually independent with equal probabilities. For $0<t<q^{-1/2}$ we consider the real logarithmic observable
\begin{equation}\label{eq:random}
 X_C(t)=\sum_{x\in|C|}-\log(1-\eps_xt^{\deg x}),
\end{equation}
where the sum is taken by increasing degree. It converges in $L^2$ and almost surely; absolute convergence of the linear term $\sum_{x\in|C|}\eps_xt^{\deg x}$ is not assumed. Fix a known constant $c$, independent of $q$, and set
\begin{equation}\label{eq:t}
 0<c<\frac1{\sqrt{12}},\qquad
 t_q=\frac c{\sqrt q},\qquad u_q=t_q^2=\frac{c^2}{q}.
\end{equation}
The observation is an independent copy of
\begin{equation}\label{eq:Y}
 Y_C=X_C(t_q)+W.
\end{equation}
The scale $t_q\asymp q^{-1/2}$ keeps the variance of the random logarithm bounded as $q$ grows. At this scale, changing the second Weil coefficient by one while keeping the first coefficient fixed changes the signal mean on the scale $q^{-2}$. The restriction on $c$ is a sufficient condition for the uniform separation estimate below; no optimality of this interval is asserted.

The noise $W$ is independent of all signs and is independently resampled along with them for each observation. Its distribution, the field size $q$, the genus, and $c$ are known. The input consists only of the observations: it does not disclose a defining equation for $C$, its closed points, or the coefficients of the random power series. An implementation which enumerates closed points could already reveal the desired information. We therefore study an abstract sampling oracle, without claiming a point-counting algorithm or charging finite-field operations for producing a sample.

For $0<\delta<1$, let $m_*(q,c,W,\delta)$ denote the least nonnegative integer $m$ for which some measurable estimator based on $m$ independent observations of~\eqref{eq:Y} returns $P_C$ with error at most $\delta$, uniformly over all genus-two curves over $\F_q$. The estimator may use arbitrary statistics of the entire sample. In particular, it need not estimate the observation mean, which may be undefined for the noise laws considered here.

For a real random variable $U$, write $\varphi_U(s)=\E[e^{isU}]$. By continuity at zero, every noise law admits a fixed frequency $\omega=\omega(W)$ such that
\begin{equation}\label{eq:frequency}
 0<\omega\le\frac12,\qquad |\varphi_W(\omega)|\ge\frac12.
\end{equation}
This frequency depends only on the known noise law and is chosen independently of $q$ and the unknown curve. Define
\[
 d_c=\frac14\min\{c^4,\,2c^2(1-12c^2)\}>0.
\]

\begin{theorem}\label{thm:main}
Fix $c$ as in~\eqref{eq:t}. Let $W$ be a fixed, known real random variable with a probability density $f$ satisfying $\sqrt f\in H^1(\R)$. Assume the independence and resampling specified above, and choose $\omega$ as in~\eqref{eq:frequency}. Put
\[
 \begin{gathered}
 J_W=\left\|(\sqrt f)'\right\|_{L^2(\R)}^2,\qquad
 K_W=128+8J_W,\\
 q_0(W)=\max\{3,\lceil K_W^{1/4}\rceil\}.
 \end{gathered}
\]
For every prime power $q\ge q_0(W)$ and $0<\delta<1/8$,
\begin{equation}\label{eq:mainbound}
 \frac{q^4}{2K_W}\log\frac1{4\delta}
 \le m_*(q,c,W,\delta)
 \le\left\lceil
 2^{12}\omega^{-2}d_c^{-2}q^4\log\frac4\delta
 \right\rceil.
\end{equation}
No moment assumption is imposed on $W$. The upper bound holds for all prime powers $q\ge2$ and every known real noise law, without any density assumption. The lower bound already holds when the admissible curves are restricted to those whose Jacobians are ordinary and absolutely simple.
\end{theorem}

Thus, for fixed $c$ and $W$ satisfying the theorem, the minimax sample complexity is $\Theta_{c,W}(q^4\log(1/\delta))$, uniformly for $0<\delta<1/8$. The constants are conservative. At $c=1/4$, one has $d_c=2^{-10}$, so the upper bound in~\eqref{eq:mainbound} becomes $\lceil2^{32}\omega^{-2}q^4\log(4/\delta)\rceil$. The finite-variance estimators in Section~\ref{sec:upper} give alternative explicit constants.

Two mechanisms produce the matching bounds. The identity
\[
 \E[X_C(t)]=\frac12\log\zeta_C(t^2)
\]
identifies a signal location parameter. Integrality of the two Weil coefficients separates its possible values by a constant times $q^{-2}$. The centered random logarithm is symmetric, and its characteristic function is positive at a sufficiently small fixed frequency. Consequently, the signal location is encoded in the phase of that characteristic function. Estimating this phase uses bounded random variables and attains the required accuracy without moments of the observation noise.

For the lower bound, we use two curves with the same first Weil coefficient and second coefficients differing by two. Their entire observation laws have squared Hellinger distance $O_W(q^{-4})$, uniformly over the stated interval of $c$. Controlling the degree-two sign block at this precision, together with the infinite closed-point tail, proves that no estimator based on the full sample can improve the order of the upper bound. The obstruction persists on the ordinary absolutely simple Jacobian locus.

The arithmetic realization uses Howe--Nart--Ritzenthaler~\cite{HNR} and Howe--Zhu~\cite{HZ}. The $b=1$ family used below already appears immediately after Theorem~6 of~\cite{HZ}; we pair it with the $b=-1$ family to obtain the required closed-point differences. The binomial identity is an elementary specialization of the Krawtchouk--Parseval framework; see Zacharovas~\cite[Section~2.1]{Zacharovas}. We retain a direct proof for the particular comparison needed here. Robust mean estimation and Hellinger testing are classical~\cite{DLLO,LeCamYang}. The contribution is the matching recovery and full-observation lower bounds for this arithmetic experiment, including their realization by ordinary absolutely simple Jacobians.

Other zeta reconstruction and distribution problems use different data. Bilu--Das--Howe~\cite{BDH} study zeta statistics and convergence in arithmetic families. Chidambaram--Keller~\cite{CK} recover Frobenius polynomials from directly supplied point counts of abelian varieties. Neither observation scheme supplies samples of the scalar random Euler channel considered here.

\section{The random logarithm and its first two moments}\label{sec:moments}
For $n\ge1$ write $N_n(C)$ for the number of closed points of degree $n$, and $T_n(C)=\#C(\F_{q^n})$. The basic relations are
\begin{equation}\label{eq:counts}
 T_n(C)=\sum_{d\mid n}dN_d(C),\qquad
 N_n(C)=\frac1n\sum_{d\mid n}\mu(n/d)T_d(C),
\end{equation}
where $\mu$ denotes the ordinary M\"obius function. For $x=q^n\ge2$, the Weil bound implies
\begin{equation}\label{eq:Weilcount}
 nN_n(C)\le T_n(C)\le q^n+1+4q^{n/2}<5q^n.
\end{equation}
Define $h_n(t)=\at(t^n)$ and let $S_n(C)$ be a sum of $N_n(C)$ independent symmetric signs, with independent sums for different $n$.

\begin{lemma}\label{lem:decomp}
For every $0<t<q^{-1/2}$, the series defining $X_C(t)$ converges in $L^2$ and almost surely, and
\begin{equation}\label{eq:decomp}
 X_C(t)=\frac12\log\zeta_C(t^2)+\sum_{n\ge1}h_n(t)S_n(C).
\end{equation}
Consequently,
\begin{equation}\label{eq:moment}
 \rev{\E[X_C(t)]}=\tfrac12\log\zeta_C(t^2),\qquad
 \Var(X_C(t))=\sum_{n\ge1}N_n(C)h_n(t)^2.
\end{equation}
\begin{revision}
In particular, $\Var(X_C(t_q))<2$ uniformly over $q\ge2$, all $c$ in~\eqref{eq:t}, and all such curves.
\end{revision}
\end{lemma}
\begin{proof}
For $0<a<1$ and $\eps\in\{-1,1\}$, we have
the exact real identity
\[
 -\log(1-\eps a)
 =-\tfrac12\log(1-a^2)+\eps\at(a).
\]
Applying this identity to~\eqref{eq:random} and
grouping the terms by degree, we obtain
\[
 X_C(t)
 =-\frac12\sum_{x\in|C|}
   \log(1-t^{2\deg x})
 +\sum_{n\ge1}h_n(t)S_n(C),
\]
where the equality is initially understood
at the level of partial sums.
Since $qt^2<1$, the Euler product for
$\zeta_C(t^2)$ converges absolutely, and hence
\[
 -\frac12\sum_{x\in|C|}
 \log(1-t^{2\deg x})
 =\frac12\log\zeta_C(t^2).
\]

The random variables $S_n(C)$ are independent
across degrees and satisfy
\[
 \rev{\E[S_n(C)]}=0,\qquad
 \Var(S_n(C))=N_n(C).
\]
Using $h_n(t)\le t^n/(1-t^2)$
and~\eqref{eq:Weilcount}, we obtain
\[
 \sum_{n\ge1}N_n(C)h_n(t)^2
 \le\frac5{(1-t^2)^2}
 \sum_{n\ge1}\frac{(qt^2)^n}{n}
 <\infty.
\]
Consequently, the partial sums
\[
 M_N=\sum_{n=1}^N h_n(t)S_n(C)
\]
form an $L^2$-bounded martingale.
By the $L^2$ martingale convergence theorem,
$M_N$ converges both in $L^2$ and almost
surely. This proves~\eqref{eq:decomp}.

Since the terms $h_n(t)S_n(C)$ are centered
and mutually orthogonal in $L^2$,
\[
 \rev{\E[X_C(t)]}
 =\frac12\log\zeta_C(t^2),
 \qquad
 \Var(X_C(t))
 =\sum_{n\ge1}N_n(C)h_n(t)^2.
\]
Thus~\eqref{eq:moment} follows.

By \eqref{eq:Weilcount},
$t_q=cq^{-1/2}$, $h_n(t_q)\le2t_q^n$, and $0<c<1/\sqrt{12}$,
we obtain
\[
\begin{aligned}
 \Var(X_C(t_q))
 &\le20\sum_{n\ge1}\frac{c^{2n}}n\\
&=20\log\frac1{1-c^2}\\
&<20\log\frac{12}{11}\\
&<2.
\end{aligned}
\]
\end{proof}

The symmetry in the decomposition makes the signal mean accessible even when the observation has no expectation. We record the precise bounds used by the estimator.
\begin{corollary}\label{cor:phase}
Put $\mu_C=\E[X_C(t_q)]$ and $Z_C=X_C(t_q)-\mu_C$. Then $Z_C\stackrel{d}{=}-Z_C$, where
$\stackrel{d}{=}$ denotes equality in distribution,
\[
 \E[Z_C^2]<2,\qquad
 0\le\mu_C\le\frac52\log\frac1{1-c^2}<\frac14.
\]
For every real $s$ with $|s|\le1/2$,
\begin{equation}\label{eq:signalphase}
 \varphi_{X_C(t_q)}(s)=e^{is\mu_C}A_C(s),\qquad
 A_C(s)=\E[\cos(sZ_C)]\ge1-s^2\ge\frac34.
\end{equation}
For fixed $q,c$, the law of $X_C(t_q)$ depends on $C$ only through $P_C$.
\end{corollary}
\begin{proof}
By~\eqref{eq:decomp}, we have
$Z_C=\sum_{n\ge1}h_n(t_q)S_n(C)$.
Since the random variables $S_n(C)$ are
independent and symmetric, every finite
partial sum of this series has a symmetric
distribution. By Lemma~\ref{lem:decomp},
these partial sums converge in $L^2$,
and hence also in distribution.
Therefore $Z_C\stackrel{d}{=}-Z_C$.
Moreover, since $\E[Z_C]=0$,
Lemma~\ref{lem:decomp} gives
$\E[Z_C^2]=\Var(X_C(t_q))<2$.
By the expectation formula
in~\eqref{eq:moment}, the definition
$u_q=c^2/q$, and the Weil bound
in~\eqref{eq:Weilcount}, we obtain
\[
 0\le\mu_C
 =\frac12\sum_{n\ge1}
   \frac{T_n(C)u_q^n}{n}\le\frac52\sum_{n\ge1}
   \frac{c^{2n}}n\frac52\log\frac1{1-c^2}<\frac52\log\frac{12}{11}
 <\frac14.
\]
Since $Z_C$ has a symmetric distribution
and $\sin$ is odd, we have
$\E[\sin(sZ_C)]=0$ for every real $s$.
Consequently,
\[
\begin{aligned}
 \varphi_{X_C(t_q)}(s)
 &=\E[e^{is(\mu_C+Z_C)}]\\
 &=e^{is\mu_C}\E[e^{isZ_C}]\\
 &=e^{is\mu_C}\E[\cos(sZ_C)]\\
 &=e^{is\mu_C}A_C(s).
\end{aligned}
\]
Furthermore, the inequality
$\cos x\ge1-x^2/2$ for $x\in\R$
implies
\[
\begin{aligned}
 A_C(s)
 &\ge1-\frac{s^2}{2}\E[Z_C^2]\\
 &\ge1-s^2\\
&\ge\frac34
\end{aligned}
\]
whenever $|s|\le1/2$.
This proves~\eqref{eq:signalphase}.

Finally, for fixed $q$ and $c$,
the zeta numerator $P_C$ determines
$\zeta_C$ and hence all point counts
$T_n(C)$.
By M\"obius inversion in~\eqref{eq:counts},
it also determines all closed-point
counts $N_n(C)$.
Since each $S_n(C)$ is a sum of
$N_n(C)$ independent symmetric signs,
and these sums are independent across
degrees, their joint distribution is
determined by $(N_n(C))_{n\ge1}$.
Thus~\eqref{eq:decomp} shows that
the law of $X_C(t_q)$ depends on $C$
only through $P_C$.
\end{proof}

\section{Arithmetic separation and recovery estimators}\label{sec:upper}

\subsection{Arithmetic separation}
Let $\mathscr W_q$ be the finite set of polynomials of the form~\eqref{eq:poly} with integral coefficients and all reciprocal roots of absolute value $\sqrt q$. Some elements of $\mathscr W_q$ need not be realized by curves; allowing them as candidates only makes the estimator less dependent on a realization classification.

\begin{lemma}\label{lem:separation}
For distinct $P,Q\in\mathscr W_q$,
put $u=u_q=c^2/q$. Then
\[
 \frac12|\log P(u)-\log Q(u)|\ge d_cq^{-2}=:\Delta_{q,c}.
\]
In particular, different genus-two zeta numerators have signal means separated by at least $\Delta_{q,c}$.
\end{lemma}
\begin{proof}
Write
$P(T)=1+aT+bT^2+qaT^3+q^2T^4$
and
$Q(T)=1+a'T+b'T^2+qa'T^3+q^2T^4$.
Put $\Delta a=a-a'$ and $\Delta b=b-b'$.
Then
\begin{equation}\label{eq:diffP}
 P(u)-Q(u)
 =\Delta a\,u(1+qu^2)+\Delta b\,u^2.
\end{equation}

We first estimate the coefficients.
Let $\alpha_1,\ldots,\alpha_4$ be the reciprocal
roots of $P$. By the Weil bounds,
$|\alpha_i|=\sqrt q$ for every $i$.
Since $b=\sum_{i<j}\alpha_i\alpha_j$,
we have $|b|\le6q$.
Similarly, $|b'|\le6q$, and hence
$|\Delta b|\le12q$.

We distinguish two cases.

\smallskip
\noindent
\emph{Case 1: $\Delta a=0$.}
Since $P\ne Q$, we have $\Delta b\ne0$.
As $\Delta b$ is an integer, $|\Delta b|\ge1$.
Thus, by~\eqref{eq:diffP} and $u=c^2/q$,
we obtain
$|P(u)-Q(u)|=|\Delta b|u^2
\ge u^2=c^4q^{-2}$.

\smallskip
\noindent
\emph{Case 2: $\Delta a\ne0$.}
Since $\Delta a$ is a nonzero integer,
$|\Delta a|\ge1$.
Using the reverse triangle inequality
and $|\Delta b|\le12q$, we obtain
\[
\begin{aligned}
 |P(u)-Q(u)|
 &\ge |\Delta a|u(1+qu^2)-|\Delta b|u^2\\
 &\ge u(1+qu^2)-12qu^2\\
 &=\frac{c^2}{q}
   \left(1-12c^2+\frac{c^4}{q}\right)\\
 &\ge c^2(1-12c^2)q^{-1}\\
 &\ge 2c^2(1-12c^2)q^{-2}.
\end{aligned}
\]
Here we used $1-12c^2>0$ and $q\ge2$.

Combining the two cases, we conclude that
\[
 |P(u)-Q(u)|
 \ge \min\{c^4,\,2c^2(1-12c^2)\}q^{-2}.
\]

It remains to estimate the difference of
the logarithms.
Write $P(T)=\prod_{i=1}^4(1-\alpha_iT)$,
where $|\alpha_i|=\sqrt q$.
Since $\sqrt q\,u=c^2/\sqrt q<1$,
we have $P(t)\ne0$ for every $0\le t\le u$.
As $P$ has real coefficients and $P(0)=1$,
continuity implies $P(u)>0$.
Moreover, by the Weil bounds,
\[
 P(u)\le\prod_{i=1}^4(1+|\alpha_i|u)=(1+\sqrt q\,u)^4=\left(1+\frac{c^2}{\sqrt q}\right)^4\le\left(1+\frac1{12\sqrt2}\right)^4<2.
\]
The same argument applies to $Q$.
Consequently, $0<P(u),Q(u)<2$.
By the mean value theorem applied to
$f(x)=\log x$, there exists a number
$\xi$ between $P(u)$ and $Q(u)$ such that
$\log P(u)-\log Q(u)=(P(u)-Q(u))/\xi$.
Since $0<\xi<2$, it follows that
\[
\begin{aligned}
 \frac12|\log P(u)-\log Q(u)|
 &=\frac{|P(u)-Q(u)|}{2\xi}\\
 &\ge\frac14|P(u)-Q(u)|\\
 &\ge\frac14
   \min\{c^4,\,2c^2(1-12c^2)\}q^{-2}\\
 &=d_cq^{-2}.
\end{aligned}
\]

Finally, since
$\zeta_C(T)=P_C(T)/((1-T)(1-qT))$,
the signal means corresponding to $P$ and $Q$
differ by
$\frac12|\log P(u)-\log Q(u)|$.
This proves the asserted separation.
\end{proof}

\subsection{Recovery from a characteristic-function phase}
We first recover the signal location from a bounded statistic. This construction needs only the known noise law and its characteristic function near the origin.

\begin{proposition}\label{prop:phaseupper}
Let $W$ have any known real probability law, with the independence and resampling specified in the model. Choose $\omega$ as in~\eqref{eq:frequency}. For every prime power $q\ge2$ and $0<\delta<1$, there is an estimator of $P_C$ using at most
\begin{equation}\label{eq:phasebudget}
 \left\lceil2^{12}\omega^{-2}d_c^{-2}q^4\log\frac4\delta\right\rceil
\end{equation}
independent observations and having error at most $\delta$, uniformly over all genus-two curves over $\F_q$.
\end{proposition}
\begin{proof}
Let $\Arg$ denote the principal argument in $(-\pi,\pi]$, and put
$\theta_W=\Arg\varphi_W(\omega)$. 
Recall that for a real random variable $U$,
its characteristic function is defined by
$\varphi_U(s)=\E[e^{isU}]$.
Since $Y_C=X_C(t_q)+W$ and
$X_C(t_q)$ is independent of $W$,
we have
$\varphi_{Y_C}(\omega)
=\varphi_{X_C(t_q)}(\omega)\varphi_W(\omega)$.
By Corollary~\ref{cor:phase},
$\varphi_{X_C(t_q)}(\omega)
=e^{i\omega\mu_C}A_C(\omega)$.
Since $\theta_W=\Arg\varphi_W(\omega)$
and $\varphi_W(\omega)\ne0$, we have
$\varphi_W(\omega)
=|\varphi_W(\omega)|e^{i\theta_W}$.
Consequently,
\[
\begin{aligned}
 z_C
 &:=e^{-i\theta_W}\varphi_{Y_C}(\omega)\\
 &=e^{-i\theta_W}e^{i\omega\mu_C}
   A_C(\omega)|\varphi_W(\omega)|
   e^{i\theta_W}\\
 &=|\varphi_W(\omega)|A_C(\omega)
   e^{i\omega\mu_C}\\
 &=r_Ce^{i\omega\mu_C},
\end{aligned}
\]
where
$r_C=|\varphi_W(\omega)|A_C(\omega)
\ge(1/2)(3/4)=3/8$.
Moreover, since $0\le\mu_C<1/4$
by Corollary~\ref{cor:phase}
and $0<\omega\le1/2$
by~\eqref{eq:frequency},
we have $0\le\omega\mu_C<1/8$.

Define
\[
 \widehat z=\frac1m\sum_{j=1}^m
 e^{i(\omega Y_j-\theta_W)}.
\]
The real and imaginary parts of
$e^{i(\omega Y_j-\theta_W)}$
are $\cos(\omega Y_j-\theta_W)$ and
$\sin(\omega Y_j-\theta_W)$, respectively,
and hence both lie in $[-1,1]$.
Since $\E[\widehat z]=z_C$,
Hoeffding's inequality gives, for every
$t>0$,
\[
\begin{aligned}
 \Prob\{|\operatorname{Re}(\widehat z-z_C)|>t\}
 &\le 2e^{-mt^2/2},\\
 \Prob\{|\operatorname{Im}(\widehat z-z_C)|>t\}
 &\le 2e^{-mt^2/2}.
\end{aligned}
\]
If $|\widehat z-z_C|>r$, then either
$|\operatorname{Re}(\widehat z-z_C)|>r/\sqrt2$
or
$|\operatorname{Im}(\widehat z-z_C)|>r/\sqrt2$.
Thus, by the union bound with
$t=r/\sqrt2$, we obtain
\begin{equation}\label{eq:ecfconcentration}
 \Prob\{|\widehat z-z_C|>r\}
 \le 4\exp(-mr^2/4).
\end{equation}

Set $r=\omega\Delta_{q,c}/32$.
Since $d_c\le c^4/4<1/576$,
$\omega\le1/2$, and $q\ge2$, we have
$r<3/16\le|z_C|/2$.
On the event $|\widehat z-z_C|\le r$,
define $v=(\widehat z-z_C)/z_C$.
Then $\widehat z=z_C(1+v)$ and
$|v|\le r/|z_C|<1/2$.
In particular,
$\operatorname{Re}(1+v)\ge1-|v|>0$,
so $1+v$ lies in the right half-plane.
Using $|\arctan x|\le|x|$ for real $x$,
we obtain
\[
\begin{aligned}
 |\Arg(1+v)|
 &=\left|\arctan
   \frac{\operatorname{Im}v}
        {1+\operatorname{Re}v}\right|\\
 &\le\frac{|\operatorname{Im}v|}
           {1+\operatorname{Re}v}\\
 &\le\frac{|v|}{1-|v|}\\
 &\le2|v|.
\end{aligned}
\]
Together with $0\le\omega\mu_C<1/8$ and $8r<1/8$, this also shows that the principal arguments do not cross a branch cut. Consequently,
\begin{equation}\label{eq:phaseerror}
 |\Arg\widehat z-\omega\mu_C|
 \le\frac{2r}{|z_C|}\le8r
 =\frac{\omega\Delta_{q,c}}4.
\end{equation}

Define
\[
 \widehat\mu=
 \begin{cases}
  \Arg\widehat z/\omega,
    &\text{if }\widehat z\ne0,\\
  0,
    &\text{if }\widehat z=0.
 \end{cases}
\]
For each $P\in\mathscr W_q$, put
\begin{equation}\label{eq:candidatelocation}
 \mu_P=\frac12\log
 \frac{P(u_q)}{(1-u_q)(1-qu_q)}.
\end{equation}
We select a polynomial $P\in\mathscr W_q$
minimizing $|\widehat\mu-\mu_P|$,
resolving ties by a fixed order.
On the event $|\widehat z-z_C|\le r$,
we have
$|\widehat z|\ge|z_C|-r>0$.
Thus~\eqref{eq:phaseerror} gives
$|\widehat\mu-\mu_C|
\le\Delta_{q,c}/4$.
By Lemma~\ref{lem:separation},
the true numerator $P_C$ is
the unique nearest candidate.

Since $r=\omega d_cq^{-2}/32$, we have
\[
 \frac4{r^2}
 =2^{12}\omega^{-2}d_c^{-2}q^4.
\]
Thus, choosing $m$ as in~\eqref{eq:phasebudget},
inequality~\eqref{eq:ecfconcentration} gives
\[
 \Prob\{|\widehat z-z_C|>r\}
 \le4\exp(-mr^2/4)
 \le\delta.
\]
Consequently, the estimator recovers
$P_C$ with probability at least $1-\delta$,
as required.
The concentration argument uses only
the bounded real and imaginary parts of
$e^{i(\omega Y_j-\theta_W)}$.
In particular, no moment assumption
on $W$ or $Y_C$ is needed.
\end{proof}

\subsection{Estimators for finite-variance noise}
When the noise has finite variance, an alternative estimator uses only block averages. It can be implemented from the noise mean and a variance bound, without evaluating a characteristic function.

\begin{proposition}\label{prop:upper}
Suppose the noise mean $m_W=\E[W]$ and a finite upper bound $v_W\ge\Var(W)$ are known. For $q\ge2$ and $0<\delta<1/8$, there is an estimator of $P_C$ based on at most
\[
 2^{10}(2+v_W)d_c^{-2}q^4\log(2/\delta)
\]
independent observations whose error is at most $\delta$, uniformly over all genus-two curves over $\F_q$.
\end{proposition}
\begin{proof}
Set $V_W=2+v_W$ and
$\epsilon_{q,c}=\Delta_{q,c}/4$.
Since $X_C(t_q)$ and $W$ are independent,
we have
$\Var(Y_C)<2+v_W=V_W$.
Let $B$ be the smallest positive odd integer
satisfying $B\ge\log_2(1/\delta)$, and set
\[
 L=\left\lceil16V_W\epsilon_{q,c}^{-2}\right\rceil
  =\left\lceil2^8V_Wd_c^{-2}q^4\right\rceil.
\]
Take $BL$ independent observations of $Y_C$
and divide them into $B$ disjoint blocks,
each containing $L$ observations.
For $1\le b\le B$, let $M_b$ denote
the sample mean of the observations
in the $b$-th block.
Define $\widehat\nu$ to be the median
of $M_1,\ldots,M_B$.
Since the observations are independent
and identically distributed, we have
$\E[M_b]=\E[Y_C]$ and
$\Var(M_b)\le V_W/L$.
By Chebyshev's inequality,
\[
\begin{aligned}
 \Prob\{|M_b-\E[Y_C]|>\epsilon_{q,c}\}&\le
 \frac{\Var(M_b)}{\epsilon_{q,c}^2}\\
 &\le
 \frac{V_W}{L\epsilon_{q,c}^2}\\
&\le\frac1{16}.
\end{aligned}
\]

We now estimate the probability that
$\widehat\nu$ differs from $\E[Y_C]$
by more than $\epsilon_{q,c}$.
Since $B$ is odd, this can happen only if
at least $k=(B+1)/2$ block means differ
from $\E[Y_C]$ by more than $\epsilon_{q,c}$.
The blocks are independent, so the
probability that any specified $k$
block means all have such errors
is at most $(1/16)^k$.
There are $\binom{B}{k}$ possible
choices of these $k$ blocks.
Thus the union bound gives
\begin{equation}\label{eq:MoM}
\begin{aligned}
 \Prob\{|\widehat\nu-\E[Y_C]|
          >\epsilon_{q,c}\}&\le
 \binom{B}{k}\left(\frac1{16}\right)^k\\
 &\le
 2^B\left(\frac1{16}\right)^{B/2}\\
 &=2^{-B}\le\delta.
\end{aligned}
\end{equation}
Recall that $\E[Y_C]=m_W+\mu_C$.
For each candidate polynomial
$P\in\mathscr W_q$, let $\mu_P$
be defined by~\eqref{eq:candidatelocation}.
We define the estimator to be a polynomial
$P\in\mathscr W_q$ that minimizes
$|\widehat\nu-(m_W+\mu_P)|$,
resolving ties by a fixed order.
Consider the event
$|\widehat\nu-\E[Y_C]|
\le\epsilon_{q,c}$.
On this event, the true numerator $P_C$
satisfies
\[
 |\widehat\nu-(m_W+\mu_C)|
 \le\frac{\Delta_{q,c}}4.
\]
For any other candidate $P\ne P_C$,
Lemma~\ref{lem:separation} and
the triangle inequality give
\[
\begin{aligned}
 |\widehat\nu-(m_W+\mu_P)|
 &\ge|\mu_P-\mu_C|
     -|\widehat\nu-(m_W+\mu_C)|\\
 &\ge\Delta_{q,c}
     -\frac{\Delta_{q,c}}4\\
 &=\frac{3\Delta_{q,c}}4.
\end{aligned}
\]
Therefore, $P_C$ is the unique candidate
closest to $\widehat\nu$.
By~\eqref{eq:MoM}, the estimator
recovers $P_C$ with probability
at least $1-\delta$.
No union bound over the candidate
polynomials is needed.

It remains to verify the sample size.
By the choice of $B$,
we have $B\le3\log(2/\delta)$.
Moreover, since $V_Wd_c^{-2}q^4\ge1$,
the definition of $L$ gives
\[
 L\le(2^8+1)V_Wd_c^{-2}q^4.
\]
Hence the total number of observations
satisfies
\[
\begin{aligned}
 LB
 &\le3(2^8+1)V_Wd_c^{-2}q^4
       \log(2/\delta)\\
 &<2^{10}V_Wd_c^{-2}q^4
       \log(2/\delta).
\end{aligned}
\]
This proves the proposition.
\end{proof}

For Gaussian noise, the centered signal and the noise are both sub-Gaussian, so the ordinary average of all observations suffices.
\begin{corollary}\label{cor:gaussmean}
For $q\ge2$, $0<\delta<1$, and $W\sim N(\mu_W,\sigma^2)$ with $\sigma>0$, the ordinary sample mean followed by nearest-candidate selection succeeds with error at most $\delta$ whenever
\[
 m\ge2^5(2+\sigma^2)d_c^{-2}q^4\log(2/\delta).
\]
\end{corollary}
\begin{proof}
Recall that $Z_C=X_C(t_q)-\mu_C$.
By~\eqref{eq:decomp}, we have
$Z_C=\sum_{n\ge1}h_n(t_q)S_n(C)$,
where the random variables $S_n(C)$
are independent, and each $S_n(C)$
is a sum of $N_n(C)$ independent
symmetric signs.
For a positive integer $N$, define
$Z_C^{(N)}=\sum_{n=1}^N h_n(t_q)S_n(C)$.
Since each symmetric sign takes the
values $1$ and $-1$ with equal
probability, its exponential moment
is $\cosh s$.
Using independence and the inequality
$\cosh x\le e^{x^2/2}$ for $x\in\R$,
we obtain, for every $s\in\R$,
\[
\begin{aligned}
 \E[e^{sZ_C^{(N)}}]
 &=\prod_{n=1}^N
   \cosh(sh_n(t_q))^{N_n(C)}\\
 &\le
 \exp\left(
 \frac{s^2}{2}
 \sum_{n=1}^N N_n(C)h_n(t_q)^2
 \right)\\
 &\le
 \exp\left(
 \frac{s^2}{2}\Var(X_C(t_q))
 \right).
\end{aligned}
\]
By Lemma~\ref{lem:decomp},
the partial sums $Z_C^{(N)}$
converge almost surely to $Z_C$.
Fatou's lemma therefore gives
\[
 \E[e^{sZ_C}]
 \le
 \exp\left(
 \frac{s^2}{2}\Var(X_C(t_q))
 \right)
 \le e^{s^2},
 \qquad s\in\R,
\]
where the last inequality follows from
$\Var(X_C(t_q))<2$.
Since $W\sim N(\mu_W,\sigma^2)$,
its centered exponential moment is
\[
 \E[e^{s(W-\mu_W)}]
 =e^{s^2\sigma^2/2}.
\]
Moreover, $Z_C$ and $W$ are independent.
Since
$Y_C-\E[Y_C]=Z_C+(W-\mu_W)$,
we obtain
\[
\begin{aligned}
 \E[e^{s(Y_C-\E[Y_C])}]
 &=\E[e^{sZ_C}]
   \E[e^{s(W-\mu_W)}]\\
 &\le
 \exp\left(
 \frac{s^2}{2}(2+\sigma^2)
 \right).
\end{aligned}
\]
Let $Y_1,\ldots,Y_m$ be independent
observations of $Y_C$, and define
their sample mean by
$\overline Y=m^{-1}\sum_{j=1}^mY_j$.
By independence, the preceding bound
implies
\[
 \E[e^{s(\overline Y-\E[Y_C])}]
 \le
 \exp\left(
 \frac{s^2(2+\sigma^2)}{2m}
 \right).
\]
For every $\epsilon>0$ and $s>0$,
Markov's inequality gives
\[
\begin{aligned}
 \Prob\{\overline Y-\E[Y_C]>\epsilon\}&\le
 e^{-s\epsilon}
 \E[e^{s(\overline Y-\E[Y_C])}]\\
 &\le
 \exp\left(
 -s\epsilon+
 \frac{s^2(2+\sigma^2)}{2m}
 \right).
\end{aligned}
\]
Choosing $s=m\epsilon/(2+\sigma^2)$
and applying the same argument
to the lower tail, we obtain
\[
 \Prob\{|\overline Y-\E[Y_C]|>\epsilon\}
 \le
 2\exp\left(
 -\frac{m\epsilon^2}{2(2+\sigma^2)}
 \right).
\]
Set $\epsilon_{q,c}=\Delta_{q,c}/4$.
Since $\Delta_{q,c}=d_cq^{-2}$, we have
$\epsilon_{q,c}^2
 =\frac{d_c^2q^{-4}}{16}$.
Hence
\[
 \Prob\{|\overline Y-\E[Y_C]|
          >\epsilon_{q,c}\}\le
 2\exp\left(
 -\frac{md_c^2q^{-4}}{2^5(2+\sigma^2)}
 \right).
\]
Thus the stated lower bound on $m$
makes this probability at most $\delta$.
For each $P\in\mathscr W_q$,
let $\mu_P$ be defined
by~\eqref{eq:candidatelocation}.
We select the candidate $P$ that
minimizes
$|\overline Y-(\mu_W+\mu_P)|$,
resolving ties by a fixed order.
On the event
$|\overline Y-\E[Y_C]|
\le\epsilon_{q,c}$,
the true numerator $P_C$ satisfies
\[
 |\overline Y-(\mu_W+\mu_C)|
 \le\frac{\Delta_{q,c}}4.
\]
For every other candidate $P\ne P_C$,
Lemma~\ref{lem:separation} and
the triangle inequality give
\[
\begin{aligned}
 |\overline Y-(\mu_W+\mu_P)|
 &\ge|\mu_P-\mu_C|
   -|\overline Y-(\mu_W+\mu_C)|\\
 &\ge\frac{3\Delta_{q,c}}4.
\end{aligned}
\]
Therefore, $P_C$ is the unique nearest
candidate on this event.
Since the event has probability at least
$1-\delta$, the estimator recovers
the true zeta numerator with error
probability at most $\delta$.
\end{proof}

\section{An arithmetic hard pair and closed-point differences}\label{sec:arithmetic}

\subsection{Realization and closed-point differences}
To prove a lower bound for all estimators, we need actual curves whose different zeta numerators generate similar observation laws. We choose the first coefficient to agree, so that the degree-one block is shared. A difference of two in the second coefficient gives degree-two sign counts on the same parity lattice. These choices lead to the pair
\begin{equation}\label{eq:hardpair}
 P_\pm(T)=1+T\pm T^2+qT^3+q^2T^4.
\end{equation}

\begin{lemma}\label{lem:hardpair}
For every prime power $q\ge3$ there are smooth, projective, geometrically connected genus-two curves $C_+,C_-$ over $\F_q$ with zeta numerators $P_+,P_-$, respectively, and both their Jacobians are ordinary and absolutely simple.
\end{lemma}
\begin{proof}
Write $b\in\{1,-1\}$ and consider the reciprocal characteristic polynomial
\[
 F_b(x)=x^4+x^3+bx^2+qx+q^2
 =x^2\bigl((x+q/x)^2+(x+q/x)+(b-2q)\bigr).
\]
Putting $s=x+q/x$, the equation $F_b(x)=0$
reduces to $s^2+s+(b-2q)=0$.
The two roots of this quadratic equation are
\[
 s_\pm=\frac{-1\pm\sqrt{8q+1-4b}}2.
\]
Consequently,
$F_b(x)=(x^2-s_+x+q)(x^2-s_-x+q)$.
Since $q\ge3$ and $b\in\{1,-1\}$, we have
$8q+1-4b>0$ and
$(1+\sqrt{8q+5})/2<2\sqrt q$.
Hence $s_+$ and $s_-$ are distinct real numbers
satisfying $-2\sqrt q<s_\pm<2\sqrt q$.
Consequently, the quadratic factors
$x^2-s_\pm x+q$ have negative discriminants
and thus have nonreal conjugate roots whose
product is $q$.
It follows that all four roots of $F_b$
have absolute value $\sqrt q$.
Then $F_b$ is a $q$-Weil polynomial.
Since $F_b$ is a $q$-Weil polynomial and $p\nmid b$,
\cite[Appendix, Step~2]{HNR} shows that $F_b$
is the Frobenius characteristic polynomial of an
ordinary abelian surface over $\F_q$.
Moreover, since $D=8q+1-4b\equiv5\pmod8$ is not a square, the corresponding isogeny class is simple.
By~\cite[Theorem~1.2(2) and Table~1.2]{HNR},
both isogeny classes contain genus-two Jacobians.
By~\cite[Theorem~6]{HZ},  both isogeny
classes are absolutely simple.
\end{proof}

The upper bound uses the separation of the signal locations determined by the two Weil coefficients. A lower bound requires more: we must compare the full laws of the independent sign blocks in every degree, including their infinite tail.

\begin{lemma}\label{lem:low}
For the curves in Lemma~\ref{lem:hardpair}, let $N_n^\pm=N_n(C_\pm)$ and $\Delta N_n=N_n^+-N_n^-$. Then
\begin{align}\label{eq:lowcounts}
 N_1^+=N_1^-&=q+2,\\
 N_2^+=\frac{q^2-q}{2},\quad& N_2^-=\frac{q^2-q-4}{2},\\
 N_3^+=\frac{q^3+2q-3}{3},\quad& N_3^-=N_3^++2,
\end{align}
and
\begin{equation}\label{eq:delta}
 \Delta N_1=0,\quad\Delta N_2=2,\quad\Delta N_3=-2,
 \quad\Delta N_4=1,\quad\Delta N_5=-2(q+1).
\end{equation}
\end{lemma}
\begin{proof}
Let $p_n(b)$ be the sum of the $n$th powers
of the roots of $F_b$. By Newton's identities,
the first five power sums are
\begin{align*}
 p_1(b)&=-1,
 &p_2(b)&=1-2b,\\
 p_3(b)&=-1+3b-3q,
 &p_4(b)&=1-4b+2b^2+4q-4q^2,\\
 p_5(b)&=-1+5b-5b^2-5q+5bq+5q^2.
\end{align*}
By the Weil point-counting formula
(see~\cite[Chapter~5]{Stichtenoth}),
$T_n(b)=q^n+1-p_n(b)$. We obtain
$T_1(b)=q+2$, $T_2(b)=q^2+2b$, and
$T_3(b)=q^3+2+3q-3b$.
By M\"obius inversion in~\eqref{eq:counts},
we have $N_1(b)=T_1(b)$,
$N_2(b)=(T_2(b)-T_1(b))/2$, and
$N_3(b)=(T_3(b)-T_1(b))/3$.
Consequently,
\[
 N_1(b)=q+2,\qquad
 N_2(b)=\frac{q^2-q+2b-2}{2},\qquad
 N_3(b)=\frac{q^3+2q-3b}{3}.
\]
Substituting $b=1$ and $b=-1$ gives
the formulas in~\eqref{eq:lowcounts}.

Newton's identities also give
$\Delta T_4=8$ and $\Delta T_5=-10(q+1)$.
Together with
$\Delta T_1=0$, $\Delta T_2=4$, and
$\Delta T_3=-6$, M\"obius inversion yields
\[
\Delta N_1=0,\quad
\Delta N_2=2,\quad
\Delta N_3=-2,\quad 
\Delta N_4=\frac{8-4}{4}=1,\quad
\Delta N_5=-\frac{10(q+1)}5=-2(q+1).
\]
This proves~\eqref{eq:delta}.
\end{proof}

\begin{lemma}\label{lem:diff}
For $q\ge3$ and every $n\ge2$,
\[
 |\Delta N_n|\le2\binom{n+2}{4}q^{n/2-1}.
\]
\end{lemma}
\begin{proof}
For each real number $b\in[-1,1]$, define
$P_b(T)=1+T+bT^2+qT^3+q^2T^4$.
Then $P_{-1}(T)=P_-(T)$ and
$P_1(T)=P_+(T)$.
As in the proof of Lemma~\ref{lem:hardpair},
the polynomial $F_b(x)=x^4P_b(1/x)$
factors as $(x^2-s_+x+q)(x^2-s_-x+q)$,
where $s_\pm=(-1\pm\sqrt{8q+1-4b})/2$.
Since $q\ge3$ and $b\in[-1,1]$, we have
$|s_\pm|<2\sqrt q$.
Thus both quadratic factors have negative
discriminants, and their roots occur in
complex-conjugate pairs with product $q$.
It follows that all four roots of $F_b$
have absolute value $\sqrt q$.
We write $[T^j]f(T)$ for the coefficient
of $T^j$ in a power series $f(T)$.
Since
$P_b(T)^{-1}=\prod_{i=1}^4(1-\alpha_iT)^{-1}$, we have
\[
 [T^j]P_b(T)^{-1}
 =\sum_{\substack{k_1+\cdots+k_4=j\\k_i\ge0}}
   \alpha_1^{k_1}\cdots\alpha_4^{k_4}.
\]
There are $\binom{j+3}{3}$ terms in this sum,
and each has absolute value $q^{j/2}$.
Hence
\begin{equation}\label{eq:coeffbound}
 |[T^j]P_b(T)^{-1}|
 \le\binom{j+3}{3}q^{j/2}.
\end{equation}
By the Weil point-counting formula,
$\Delta T_n=p_n(-1)-p_n(1)$.
Since $\partial_b P_b(T)=T^2$ and
$\log P_b(T)=-\sum_{n\ge1}p_n(b)T^n/n$,
integration with respect to $b$ gives
\[
\int_{-1}^1\frac{T^2}{P_b(T)}\,db=\log\frac{P_+(T)}{P_-(T)}=\sum_{n\ge1}\frac{\Delta T_n}{n}T^n.
\]
The power series converge uniformly in
$b\in[-1,1]$ on every closed disk
$|T|\le r<q^{-1/2}$, which justifies
termwise integration.
Comparing coefficients of $T^n$ for
$n\ge2$, we obtain
\[
 \frac{\Delta T_n}{n}
 =\int_{-1}^1[T^{n-2}]P_b(T)^{-1}\,db.
\]
By~\eqref{eq:coeffbound}, this gives
$|\Delta T_n|
\le2n\binom{n+1}{3}q^{n/2-1}$.
Since $\Delta T_1=0$, M\"obius inversion
in~\eqref{eq:counts} yields
\begin{align*}
|\Delta N_n|
&\le\frac1n\sum_{d\mid n}|\Delta T_d|\\
&\le\frac2n\sum_{\substack{d\mid n\\d\ge2}}
 d\binom{d+1}{3}q^{d/2-1}\\
&\le2q^{n/2-1}
 \sum_{d=2}^n\binom{d+1}{3}\\
&=2\binom{n+2}{4}q^{n/2-1}.
\end{align*}
Here we used $d\le n$, $q\ge3$, and the
binomial identity
$\sum_{d=2}^n\binom{d+1}{3}=\binom{n+2}{4}$.
\end{proof}

\subsection{Uniform tail and location bounds}
The exact differences in degrees four and five
are needed for the sharp tail order.
The general bound alone would lose a power of
$q$ in degree four, whereas treating these two
degrees separately preserves the $q^{-4}$
squared-coupling scale.

\begin{lemma}\label{lem:tail}
Let $h_n=h_n(t_q)$ and $\mu_\pm=\E[X_{C_\pm}(t_q)]$. For every $q\ge3$ and every $c$ in~\eqref{eq:t},
\[
 \sum_{n\ge4}|\Delta N_n|h_n^2<q^{-4},
\]
and
\[
 0<\mu_+-\mu_-<u_q^2=c^4q^{-2}.
\]
\end{lemma}
\begin{proof}
Since $t_q=c/\sqrt q<1/6$ and
$\operatorname{arctanh}(x)\le x/(1-x^2)<2x$
for $0<x<1/6$, we have
$h_n=\operatorname{arctanh}(t_q^n)<2t_q^n$.
Consequently,
$h_n^2\le4c^{2n}q^{-n}$ for every $n\ge1$.
By~\eqref{eq:delta}, we have
$|\Delta N_4|=1$ and
$|\Delta N_5|=2(q+1)$.
Using $h_n^2\le4c^{2n}q^{-n}$,
we obtain
\[
 |\Delta N_4|h_4^2\le4c^8q^{-4},
 \qquad
 |\Delta N_5|h_5^2
 \le8(q+1)c^{10}q^{-5}
 \le\frac{32}{3}c^{10}q^{-4},
\]
where the last inequality follows from $q\ge3$.
For $n\ge6$, Lemma~\ref{lem:diff} gives
\[
 |\Delta N_n|h_n^2
 \le8\binom{n+2}{4}c^{2n}q^{-n/2-1}
 \le8\binom{n+2}{4}c^{2n}q^{-4}.
\]
Splitting the sum into the terms $n=4$,
$n=5$, and $n\ge6$, and using
$\sum_{n\ge2}\binom{n+2}{4}x^n=x^2/(1-x)^5$,
we obtain
\begin{align*}
 q^4\sum_{n\ge4}|\Delta N_n|h_n^2
 &\le4c^8+\frac{32}{3}c^{10}
       +\frac{8c^4}{(1-c^2)^5}\\
 &<\frac4{12^4}+\frac{32}{3\cdot12^5}
       +\frac8{12^2}\left(\frac{12}{11}\right)^5<1.
\end{align*}

It remains to estimate $\mu_+-\mu_-$.
Since $\mu_\pm=\frac12\log\zeta_{C_\pm}(u_q)$
and the two zeta functions have the same
denominator, we have
$\mu_+-\mu_-=\frac12\log(P_+(u_q)/P_-(u_q))$.
Moreover,
$P_+(u_q)-P_-(u_q)=2u_q^2$, and
$P_-(u_q)=1+u_q(1-u_q)+qu_q^3+q^2u_q^4>1$
because $0<u_q<1$.
Consequently,
\[
\begin{aligned}
0<\mu_+-\mu_-
&=\frac12\log\left(
 1+\frac{2u_q^2}{P_-(u_q)}\right)\\
&<\frac{u_q^2}{P_-(u_q)}
<u_q^2=c^4q^{-2},
\end{aligned}
\]
where we used $\log(1+x)<x$ for $x>0$.
\end{proof}

\subsection{An explicit small-field example}\label{sec:example}
The hard pair need not remain at the level of Honda--Tate existence. Over $\F_3$, consider the smooth projective models of
\begin{equation}\label{eq:curves}
 C_+:y^2=x^6+x^4+x^2+x,\qquad
 C_-:y^2=x^6+x^5+2x^4+x.
\end{equation}
The two right-hand polynomials are squarefree over $\F_3$, so these are genus-two hyperelliptic curves. They have two rational points at infinity. Direct finite-field counting gives
\[
 \#C_\pm(\F_3)=5,\qquad
 \#C_+(\F_9)=11,\qquad \#C_-(\F_9)=7.
\]
\begin{revision}
Newton's identities give $p_1=-a_C$ and $p_2=a_C^2-2b_C$, where $p_n$ denotes the Frobenius power sum. Thus $a_C=T_1(C)-(q+1)$ and $b_C=(a_C^2+T_2(C)-q^2-1)/2$. Here $a_C=1$ and $b_{C_\pm}=\pm1$, recovering~\eqref{eq:hardpair}. The following enumeration in $\F_9=\F_3[z]/(z^2+1)$ verifies these counts directly. The curves and their point counts are independent of the observation scale $c$.
\end{revision}

\subsubsection*{A hand-checkable count over \texorpdfstring{$\F_9$}{F9}}
Write $\F_9=\F_3[z]/(z^2+1)$, so $z^2=2$. Every element is $a+bz$ with $a,b\in\{0,1,2\}$. Direct evaluation of the two right-hand polynomials in~\eqref{eq:curves}, followed by the square test in $\F_9$, gives the following numbers of affine points above each $x$:
\begin{center}
\begin{tabular}{c|rrrrrrrrr}
\toprule
$x$ & $0$ & $z$ & $2z$ & $1$ & $1+z$ & $1+2z$ & $2$ & $2+z$ & $2+2z$ \\ \midrule
$C_+$ & 1 & 0 & 0 & 2 & 2 & 2 & 2 & 0 & 0 \\
$C_-$ & 1 & 0 & 0 & 2 & 0 & 0 & 2 & 0 & 0 \\ \bottomrule
\end{tabular}
\end{center}
Thus the affine totals are $9$ and $5$. Adding the two $\F_9$-rational points at infinity gives $11$ and $7$. The analogous three $x$-fibres over $\F_3$ have cardinalities $(1,2,0)$ and $(1,0,2)$, and both models have five rational points after adjoining the two points at infinity. This computation is independent of the general isogeny-class argument.

\section{Hellinger comparison and the minimax lower bound}\label{sec:hell}

\subsection{Hellinger tools and binomial comparison}
All probability measures and kernels in this section are defined on standard Borel spaces. For probability measures $P,Q$ dominated by a common measure, with densities $p,q$, use the conventions
\[
 \rho(P,Q)=\int\sqrt{pq},\qquad
 \He^2(P,Q)=\int(\sqrt p-\sqrt q)^2=2-2\rho(P,Q),
\]
and $\TV(P,Q)=\frac12\int|p-q|$. The \emph{unsquared} $\He$ obeys the triangle inequality.

For $P\ll Q$, define the chi-square divergence by $\chi^2(P\Vert Q)=\int(\dd P/\dd Q-1)^2\,\dd Q$, and set $\chi^2(P\Vert Q)=+\infty$ otherwise.

\begin{lemma}\label{lem:hellbasic}
Suppose $P_z,Q_z$ are probability kernels and $\pi$ is a common mixing distribution. Then
\[
 \He^2\left(\int P_z\dd\pi(z),\int Q_z\dd\pi(z)\right)
 \le\int\He^2(P_z,Q_z)\dd\pi(z).
\]
A common Markov kernel (in particular, a common independent convolution) contracts $\He$. Moreover,
\begin{equation}\label{eq:product}
 \rho(P^{\otimes m},Q^{\otimes m})=\rho(P,Q)^m,
 \qquad \TV(P,Q)\le\sqrt{1-\rho(P,Q)^2}.
\end{equation}
\end{lemma}
\begin{proof}
First consider a Markov kernel $K$. Put $\nu=(P+Q)/2$,
$p=dP/d\nu$, and $q=dQ/d\nu$, and let $(A,B)$ have
joint law $\nu(da)K(a,db)$. The densities of $PK$ and $QK$
with respect to $\nu K$ are $\E[p(A)\mid B]$ and
$\E[q(A)\mid B]$, respectively. Conditional Cauchy--Schwarz gives
\[
 \rho(PK,QK)
 =\E\!\left[\sqrt{\E[p(A)\mid B]\E[q(A)\mid B]}\right]
 \ge\E[\sqrt{p(A)q(A)}]
 =\rho(P,Q).
\]
Thus $K$ contracts $\He$. Apply this conclusion to the projection
$(z,x)\mapsto x$ and the joint laws $\pi(dz)P_z(dx)$ and
$\pi(dz)Q_z(dx)$. With the dominating joint measure
$\pi(dz)(P_z+Q_z)(dx)$, their squared Hellinger distance is
$\int\He^2(P_z,Q_z)\dd\pi(z)$, which proves the mixture inequality.
The product identity follows from Tonelli's theorem. Finally,
\[
 \TV(P,Q)=\tfrac12\int|\sqrt p-\sqrt q|(\sqrt p+\sqrt q)
 \le\tfrac12\sqrt{(2-2\rho)(2+2\rho)}.
\]
\end{proof}

\begin{lemma}\label{lem:smooth}
Let $W$ have density $f$ with $\sqrt f\in H^1(\R)$, and put $J_W=\| (\sqrt f)'\|_2^2$. Then for every $a,b\in\R$,
\begin{equation}\label{eq:shift}
 \He^2(\Law(W+a),\Law(W+b))\le J_W(a-b)^2.
\end{equation}
For every coupling $(U,V)$ of two real random variables with finite second moments,
\begin{equation}\label{eq:coupling}
 \He^2(\Law(U+W),\Law(V+W))\le J_W\rev{\E[(U-V)^2]},
\end{equation}
where in the marginal laws the noise is independent of each respective signal.
\end{lemma}
\begin{proof}
Take $g=\sqrt f$. The one-dimensional Sobolev representative is locally absolutely continuous. By the fundamental theorem of calculus and Cauchy--Schwarz,
\[
 \|g(\cdot-a)-g(\cdot-b)\|_2^2
 \le |a-b|\int_0^{|a-b|}\|g'(\cdot-s)\|_2^2\dd s
 =J_W(a-b)^2.
\]
This is~\eqref{eq:shift}. For~\eqref{eq:coupling}, mix~\eqref{eq:shift} over the joint distribution of $(U,V)$ and use Lemma~\ref{lem:hellbasic}. 
The coupling compares output laws only. The bound requires a finite squared difference of the signals, but no moment of the noise.
\end{proof}

We next compare two sign sums without using the readout noise. The identity below is an elementary special case of the Krawtchouk--Parseval formula in~\cite[Section~2.1]{Zacharovas}. Indeed, with $M=N+2$, the common transformation $s\mapsto(s+M)/2$ identifies the two laws with $\operatorname{Bin}(M-2,1/2)+1$ and $\operatorname{Bin}(M,1/2)$. A direct likelihood-ratio calculation is particularly short in this case.

\begin{lemma}\label{lem:chi}
Let $R_N$ be the law of a sum of $N$ independent symmetric signs. For every integer $N\ge0$,
\begin{equation}\label{eq:chi}
 \chi^2(R_N\Vert R_{N+2})=\frac2{(N+1)(N+2)},\qquad
 \He^2(R_N,R_{N+2})\le\frac2{(N+1)(N+2)}.
\end{equation}
\end{lemma}
\begin{proof}
Set $M=N+2$. On the lattice supporting $R_M$ the Radon--Nikodym derivative of $R_N$ is
\[
 r(s)=\frac{p_N(s)}{p_M(s)}=\frac{M^2-s^2}{M(M-1)},
\]
including the endpoints $s=\pm M$, where $p_N(s)=0$. For $S\sim R_M$,
\[
 \rev{\E[S^2]}=M,\qquad \rev{\E[S^4]}=3M^2-2M,
 \qquad\Var(S^2)=2M(M-1).
\]
Hence
\[
 \chi^2(R_N\Vert R_M)=\rev{\E_{R_M}\!\left[\left(\frac{M-S^2}{M(M-1)}\right)^2\right]}
 =\frac2{M(M-1)}.
\]
The Hellinger bound follows by integrating $(\sqrt r-1)^2\le(r-1)^2$. 
In the reverse direction the chi-square divergence is infinite, since $R_{N+2}$ assigns positive mass to endpoints outside the support of $R_N$.\qedhere
\end{proof}

\subsection{The full-observation lower bound}

\begin{proposition}\label{prop:oneobs}
Let $c$ and $W$ satisfy the hypotheses of Theorem~\ref{thm:main}. For every prime power $q\ge3$, the two curves of Lemma~\ref{lem:hardpair} satisfy
\begin{equation}\label{eq:oneobs}
 \He^2\bigl(\Law(Y_{C_+}),\Law(Y_{C_-})\bigr)
 \le K_Wq^{-4},\qquad K_W=128+8J_W.
\end{equation}
\end{proposition}
\begin{proof}
Decompose the observations using~\eqref{eq:decomp} into the deterministic mean, independent sign sums of each degree, and $W$. Convert the law for $C_+$ into that for $C_-$ in four steps, leaving all unchanged summands as a common independent convolution at each step.

First change the mean $\mu_+$ to $\mu_-$. Lemmas~\ref{lem:smooth} and~\ref{lem:tail} bound the squared Hellinger cost by $\rev{J_Wc^8q^{-4}}$.

Second replace the degree-three block $h_3R_{N_3^+}$ by $h_3R_{N_3^-}$, with $N_3^-=N_3^++2$. The exact comparison~\eqref{eq:chi}, contraction under convolution, and $N_3^+\ge q^3/3$ give squared cost at most $18q^{-6}$. Note that multiplication by a nonzero scalar does not change an $f$-divergence.

Third replace the degree-two block $h_2R_{N_2^+}$ by $h_2R_{N_2^-}$, where $N_2^+=N_2^-+2$. Use the orientation $\chi^2(R_{N_2^-}\Vert R_{N_2^-+2})$. For $q\ge3$,
\[
 N_2^-+1=\tfrac12(q^2-q-2)\ge\tfrac{2q^2}9,
 \quad N_2^-+2=\tfrac12(q^2-q)\ge\tfrac{q^2}3.
\]
Since Hellinger distance is symmetric,~\eqref{eq:chi} gives squared cost at most $27q^{-4}$.

Finally couple all degree-$n$ blocks for $n\ge4$ by sharing $\min(N_n^+,N_n^-)$ of their signs and adding independent signs for the extra points. Do this independently in each degree. The difference $D$ between the two coupled tails is centered; by $L^2$ convergence,
\[
 \rev{\E[D^2]}=\sum_{n\ge4}|\Delta N_n|h_n^2<q^{-4}
\]
by Lemma~\ref{lem:tail}. Translation smoothing bounds the fourth squared cost by $J_Wq^{-4}$. All compared laws, including their countable products of sign blocks, are well defined and preserve their prescribed marginals.

Let $d_1,\ldots,d_4$ be the four (unsquared) Hellinger distances. Using the triangle inequality followed by Cauchy--Schwarz,
\begin{align*}
 \He^2(\Law(Y_{C_+}),\Law(Y_{C_-}))
 &\le(d_1+d_2+d_3+d_4)^2\\
 &\le4\sum_{i=1}^4 d_i^2\\
 &\le\left[4\cdot27+72q^{-2}+\rev{4(1+c^8)J_W}\right]q^{-4}\\
 &\le(128+8J_W)q^{-4}.
\end{align*}
\end{proof}

\begin{lemma}\label{lem:test}
If $\He^2(P,Q)\le\eta\le1$, every test based on $m$ independent observations has equal-prior average error at least
\[
 \tfrac14(1-\eta/2)^{2m}\ge\tfrac14e^{-2m\eta}.
\]
\end{lemma}
\begin{proof}
The Bayes risk under equal priors is $(1-\TV(P^{\otimes m},Q^{\otimes m}))/2$. From~\eqref{eq:product} it is at least
\[
 \frac{1-\sqrt{1-\rho(P,Q)^{2m}}}{2}
 \ge\frac14\rho(P,Q)^{2m}.
\]
Use $\rho(P,Q)=1-\He^2(P,Q)/2\ge1-\eta/2$ and $\log(1-z)\ge-2z$ for $0\le z\le1/2$.
\end{proof}

\begin{proof}[Proof of Theorem~\ref{thm:main}]
The upper bound is Proposition~\ref{prop:phaseupper}.
For the lower bound, the two curves provided by Lemma~\ref{lem:hardpair} are admissible even in the ordinary absolutely simple locus. A zeta-numerator estimator with maximal error $\delta$ induces a test of $P_+$ against $P_-$ with both error probabilities at most $\delta$. For $q\ge q_0(W)$,~\eqref{eq:oneobs} gives $\He^2\le K_Wq^{-4}\le1$. Lemma~\ref{lem:test} implies
\[
 \delta\ge\frac14\exp(-2K_Wm/q^4).
\]
Solving for $m$ gives the lower inequality in~\eqref{eq:mainbound}. The same upper-bound estimator applies to the restricted curve class, which completes the restricted-locus assertion.
\end{proof}

\subsection{Gaussian comparison and parity}
A direct Gaussian-mixture calculation provides an independent check on the most delicate comparison. It does not use the binomial chi-square identity.

\begin{lemma}\label{lem:gaussian}
Let $h\in\R$, let $G\sim N(0,1)$, let $\eps_1,\eps_2$ be independent symmetric signs, and put $D=h(\eps_1+\eps_2)$, independent of $G$. If $\sigma>0$, then
\[
 \chi^2\bigl(\Law(\sigma G+D)\Vert\Law(\sigma G)\bigr)
 =\frac14\left[\cosh\left(\frac{4h^2}{\sigma^2}\right)-1\right].
\]
\end{lemma}
\begin{proof}
Let $\phi_\sigma$ denote the $N(0,\sigma^2)$ density. Completing the square gives
\[
 \int_\R\frac{\phi_\sigma(x-d)\phi_\sigma(x-e)}{\phi_\sigma(x)}\dd x
 =e^{de/\sigma^2}.
\]
For independent copies $D,D'$, the left side averaged over $D,D'$ is the chi-square divergence plus one. Since $D$ takes $0$ with probability $1/2$ and each of $\pm2h$ with probability $1/4$, the average is $3/4+\cosh(4h^2/\sigma^2)/4$.
\end{proof}

\begin{proposition}\label{prop:alt}
For fixed $c$ as in~\eqref{eq:t} and $W=\sigma G$ with fixed $\sigma>0$,
$\He^2(\Law(Y_{C_+}),\Law(Y_{C_-}))=O_{c,\sigma}(q^{-4})$ follows without Lemma~\ref{lem:chi}.
\end{proposition}
\begin{proof}
Write $R_{N+2}=R_N*R_2$. For each of the degree-two and degree-three comparisons, compare
$\Law(\sigma G+h_n(\eps_1+\eps_2))$ with $\Law(\sigma G)$,
then convolve both laws with the common degree-$n$ block and all unchanged summands.
Hellinger contraction and symmetry, Lemma~\ref{lem:gaussian}, $\He^2\le\chi^2$, and $\cosh z-1=O(z^2)$ as $z\to0$ give squared costs $O_\sigma(h_2^4)=O_{c,\sigma}(q^{-4})$ and $O_\sigma(h_3^4)=O_{c,\sigma}(q^{-6})$. The mean and tail changes have squared costs $O_{c,\sigma}(q^{-4})$ by Lemmas~\ref{lem:tail} and~\ref{lem:smooth}; here $J_W=(4\sigma^2)^{-1}$. Apply the triangle inequality to the four unsquared distances.
\end{proof}

\subsubsection*{The sharp comparison and the role of parity}
The degree-two blocks in the hard pair have about $q^2/2$ signs, with an increment of two. Sharing their common signs and adding two independent signs to the larger block gives squared signal difference $2h_2(t_q)^2=O_c(q^{-2})$. Translation smoothing applied to this coupling would therefore yield only a $q^{-2}$ bound for squared Hellinger distance, insufficient for the sharp lower bound.

The binomial comparison uses the distribution of the entire sign sum to recover the extra factor $q^{-2}$. Its unsmoothed comparison requires an even increment: $R_N$ and $R_{N+1}$ have disjoint parity supports and total variation distance one. Gaussian smoothing uses a different mechanism and has no such parity restriction. For a single symmetric sign $\eps$ independent of $G$, the same Gaussian integral gives
\[
 \chi^2\bigl(\Law(\sigma G+h\eps)\Vert\Law(\sigma G)\bigr)
 =\cosh(h^2/\sigma^2)-1
 =O_\sigma(h^4).
\]
Indeed, for independent copies $D,D'$ of $h\eps$, the chi-square divergence plus one is $\E[e^{DD'/\sigma^2}]$. The fourth-order bound follows from the symmetry and centering of the perturbation, even when only one sign is added.

\section{Noise assumptions and further directions}\label{sec:scope}\label{sec:end}
Smoothness of the readout noise enters the lower-bound mechanism, not the upper bound. For the uniform law on $[-1,1]$ one has, for $|h|\le2$,
\[
 \He^2(\Law(W+h),\Law(W))=|h|,
\]
so the quadratic shift bound~\eqref{eq:shift} fails. This is a counterexample to removing Sobolev regularity from \emph{that lemma}, not to the minimax rate itself. Atomic or noiseless readouts remain outside the lower bound proved here.

The phase estimator also separates moment assumptions from noise regularity. For standard Cauchy noise,
\[
 f(w)=\frac1{\pi(1+w^2)},\qquad
 \varphi_W(s)=e^{-|s|},\qquad
 J_W=\frac1\pi\int_\R\frac{w^2}{(1+w^2)^3}\dd w=\frac18.
\]
Thus $\sqrt f\in H^1(\R)$, although neither the mean nor the variance of $W$ exists. One can take $\omega=1/2$ in~\eqref{eq:frequency}. Theorem~\ref{thm:main} applies with $K_W=129$ and $q_0(W)=4$. It recovers $\mu_C$ through a characteristic-function phase, without defining an observation mean.

For comparison, the density $f(w)=2/[\pi(1+w^2)^2]$ has mean zero, variance one, and $J_W=1/2$, but no moment-generating function near zero. Both the phase estimator and Proposition~\ref{prop:upper} apply. The latter supplies a finite-variance alternative when only the noise mean and a variance bound are used.

The experiment fixes the binary phase law and observes one real value at a known $c/\sqrt q$ on each independent draw. Access to exact power-series coefficients, several evaluation points, or a different phase alphabet changes the available information. The lower bound here makes no assertion about such experiments. Likewise, $c$ and the noise law are fixed as $q$ varies; joint limits with vanishing $c$ or changing readout noise require separate analysis.

For each fixed observation scale in~\eqref{eq:t}, the integral Weil coefficients can be recovered at sample order $q^4\log(1/\delta)$ from a single noisy scalar channel. A bounded characteristic-function statistic attains this order without moments of the readout noise. Under the Sobolev regularity condition on the noise density, the full observation laws of an arithmetic hard pair show that the order cannot be improved. The same obstruction occurs among curves with ordinary absolutely simple Jacobians.

Extensions to rough noise, observation scales varying with $q$, or multipoint data require new distributional comparisons. In higher genus, matching bounds would also require suitable arithmetic realizations of nearby Weil data.
\section*{Statements and Declarations}
\noindent\textbf{Funding.}
This work was supported by JSPS KAKENHI Grant number 23K12507.

\medskip
\noindent\textbf{Competing interests.}
The authors declare that they have no competing interests.

\medskip
\noindent\textbf{Author contributions.}
All authors contributed to the conception of the study, the mathematical
analysis, and the preparation of the manuscript. All authors read and approved the final manuscript.

\medskip
\noindent\textbf{Data availability.}
No datasets were generated or analysed during the current study.

\end{document}